\documentclass{article}
\usepackage{amsmath,amsfonts,amssymb,amsthm,latexsym,graphicx}

\newif\ifFULL

\ifFULL
  \usepackage{color}

  \newcommand{\bluebegin}{\begingroup\color{blue}}
  \newcommand{\blueend}{\endgroup}

\fi

\allowdisplaybreaks[3]

\gdef\enddots{$\mathinner{\ldotp\ldotp\ldotp\ldotp}$\spacefactor=3000}

  \renewenvironment{thebibliography}[1]{%
    \begin{oldthebibliography}{#1}%
      \setlength{\parskip}{0.5ex}%
      \setlength{\itemsep}{0.5ex}%
  }%
  {%
    \end{oldthebibliography}%
  }

\newcommand{\bbbr}{\mathbb{R}}		

\newcommand{\const}{\mathrm{const}}		

\newcommand{\st}{\mathrel{:}}
\newcommand{\dd}{d}

\newcommand{\para}{\mathop{\parallel}\nolimits}
\newcommand{\upperLambda}{\overline{\Lambda}}
\newcommand{\lowerLambda}{\underline{\Lambda}}
\newcommand{\upperD}{\overline{D}}
\newcommand{\lowerD}{\underline{D}}

\theoremstyle{plain}
\newtheorem{lemma}{Lemma}
\newtheorem{theorem}{Theorem}
\newtheorem{proposition}{Proposition}
\theoremstyle{remark}
\newtheorem*{remark}{Remark}

\newlength{\IndentI}
\newlength{\IndentII}
\newlength{\IndentIII}
\newlength{\IndentIV}
\newlength{\WidthI}
\newlength{\WidthII}
\newlength{\WidthIII}
\newlength{\WidthIV}
\title{Jeffreys's law for general games of prediction:\\
  in search of a theory}
\author{A.~P.~Dawid and V.~G.~Vovk}

\begin{document}
\maketitle
\begin{abstract}
  We are interested in the following version of Jeffreys's law:
  if two predictors are predicting the same sequence of events
  and either is doing a satisfactory job,
  they will make similar predictions in the long run.
  We give a classification of instances of Jeffreys's law,
  illustrated with examples.
\end{abstract}

\section{Introduction}

In this paper we are interested in games of prediction
for which Jeffreys's law,
as stated in the abstract,
holds.
Specific true instances of Jeffreys's law will be referred to as
\emph{Jeffreys theorems}.

In Section~\ref{sec:taxonomy} we define several popular games of prediction
and state Jeffreys theorems for the absolute-loss, square-loss,
and bounded square-loss games.
These results serve as illustrations for our taxonomy of Jeffreys theorems;
namely,
we distinguish between Jeffreys theorems of level 1 (weakest),
level 2 (intermediate), and level 3 (strongest).

In Section~\ref{sec:reductions} we show that
in the case of so-called perfectly mixable games
there is no difference between the three levels of Jeffreys theorems.
Perfectly mixable games include, in particular,
log-loss games and the bounded square-loss game.

In the next section, Section~\ref{sec:level-2},
we state level 2 Jeffreys theorems,
which cover the log-loss and square-loss games
(not necessarily bounded).
In combination with the results of Section~\ref{sec:reductions}
this provides us with examples of level 3 Jeffreys theorems.
Some of the results in Section~\ref{sec:level-2}
are explicit inequalities,
not just statements of convergence.

The simple method of Section~\ref{sec:level-2}
does not work for the absolute-loss game.
In Section~\ref{sec:level-1} we will see that it is still possible
to prove a Jeffreys theorem for this game,
albeit only a level 1 one.

Perhaps the first instance of Jeffreys's law was proved
by Blackwell and Dubins \cite{blackwell/dubins:1962};
a pointwise version of their result was established in \cite{dawid:1985AS}.
Results similar to ours
but stated in terms of the algorithmic theory of randomness
were earlier obtained in \cite{vovk:1987} (developing \cite{kabanov/etal:1977})
and \cite{fujiwara:2008} in the case of the log-loss game,
and in \cite{vovk:2001brier}
(in essence developing \cite{skouras/dawid:1998})
in the case of the bounded square-loss game.

\section{Taxonomy and examples of Jeffreys theorems}
\label{sec:taxonomy}

A \emph{game of prediction} is a triple $(\Omega,\Gamma,\ell)$,
where $\Omega$ and $\Gamma$ are arbitrary sets,
called the \emph{outcome space} and \emph{prediction space},
respectively,
and $\ell:\Omega\times\Gamma\to\bbbr$ is called the \emph{loss function}.
The game is played according to the following perfect-information protocol.

\medskip

\noindent
\textsc{Competitive prediction protocol}

\noindent
\textbf{Players:} Nature, Predictor 1, Predictor 2, Sceptic

\noindent
\textbf{Protocol:}

\parshape=5
\IndentI   \WidthI
\IndentII  \WidthII
\IndentII  \WidthII
\IndentII  \WidthII
\IndentI   \WidthI
\noindent
FOR $n=1,2,\ldots$:\\
  Predictor 1 and Predictor 2 announce $\gamma_{n}^{[1]}\in\Gamma$
    and $\gamma_{n}^{[2]}\in\Gamma$.\\
  Sceptic announces $\tilde\gamma_n\in\Gamma$.\\
  Nature announces $\omega_n\in\Omega$.\\
END FOR

\medskip

\noindent
Three of the players,
two Predictors and one Sceptic,
are trying to predict the outcome $\omega_n$ to be announced by Nature.
Sceptic is just like another Predictor,
but he will be playing a special role in our story.
At step $n$, Predictor 1 and Predictor 2
issue predictions $\gamma_n^{[1]}$ and $\gamma_n^{[2]}$,
respectively.
The Predictors can consult each other when making the predictions,
and the pair $(\gamma_n^{[1]},\gamma_n^{[2]})$
can be regarded as their joint prediction.
After the two Predictors have announced,
Sceptic issues his own prediction $\tilde\gamma_n$.
Then Nature produces $\omega_n$.
Let $L^{[k]}_N := \sum_{n=1}^N \ell(\omega_n,\gamma^{[k]}_n)$
be the cumulative loss to time $N$ of Predictor $k$, $k=1,2$,
and similarly $\tilde L_N$ for Sceptic.

The \emph{absolute-loss game} is $(\bbbr,\bbbr,\ell)$
where $\ell(\omega,\gamma):=\lvert\omega-\gamma\rvert$.
The next proposition states our first Jeffreys theorem.

\begin{proposition}\label{prop:absolute}
  Sceptic has a strategy in the absolute-loss game
  that guarantees
  \begin{equation}\label{eq:level-1}
    \lim_{N\to\infty}
    \max
    \left(
      \frac{1}
      {
        \left|
          \gamma^{[1]}_N-\gamma^{[2]}_N
        \right|
      },
      L^{[1]}_N
      -
      \tilde L_N,
      L^{[2]}_N
      -
      \tilde L_N
    \right)
    =
    \infty.
  \end{equation}
\end{proposition}
\noindent
As usual, we set $1/0:=\infty$ in (\ref{eq:level-1}).
For the proof of Proposition~\ref{prop:absolute},
see Section~\ref{sec:level-1}.

We call (\ref{eq:level-1}),
perhaps with $\lvert\gamma^{[1]}_N-\gamma^{[2]}_N\rvert$
replaced by a different distance,
a \emph{level 1 Jeffreys theorem}.
It says that for a sufficiently distant outcome $\omega_N$, $N\gg1$,
at least one of the following three things happen:
the two Predictors' predictions $\gamma_N^{[1]}$ and $\gamma_N^{[2]}$
are close to each other;
Sceptic greatly outperforms Predictor 1 by time $N$;
Sceptic greatly outperforms Predictor 2 by time $N$.
The weakness of this statement
is that no ``stabilization''
is guaranteed along a given infinite sequence of outcomes $\omega_1\omega_2\ldots$:
it is possible that each one of the three terms of the disjunction
will be violated infinitely often.

A stronger Jeffreys theorem,
which we call a \emph{level 2 Jeffreys theorem},
would say that
\begin{equation}\label{eq:level-2}
  \lim_{N\to\infty}
  \left|
    \gamma^{[1]}_N-\gamma^{[2]}_N
  \right|
  =
  0
  \text{ or }
  \lim_{N\to\infty}
  \max
  \left(
    L^{[1]}_N
    -
    \tilde L_N,
    L^{[2]}_N
    -
    \tilde L_N
  \right)
  =
  \infty.
\end{equation}
An even stronger statement,
which we call a \emph{level 3 Jeffreys theorem},
would be
\begin{equation}\label{eq:level-3}
  \lim_{N\to\infty}
  \left|
    \gamma^{[1]}_N-\gamma^{[2]}_N
  \right|
  =
  0
  \text{ or }
  \lim_{N\to\infty}
  \left(
    L^{[1]}_N
    -
    \tilde L_N
  \right)
  =
  \infty
  \text{ or }
  \lim_{N\to\infty}
  \left(
    L^{[2]}_N
    -
    \tilde L_N
  \right)
  =
  \infty.
\end{equation}

The following two propositions give examples of level 2 and level 3 Jeffreys theorems.
The \emph{square-loss game} is $(\bbbr,\bbbr,\ell)$
where $\ell(\omega,\gamma):=(\omega-\gamma)^2$.
\begin{proposition}\label{prop:square}
  Sceptic has a strategy in the square-loss game
  that guarantees~(\ref{eq:level-2}).
\end{proposition}
The \emph{bounded square-loss game} is $([0,1],[0,1],\ell)$
where $\ell(\omega,\gamma):=(\omega-\gamma)^2$.
(We fix specific bounds, 0 and 1,
for outcomes and predictions,
but our results generalize in a straightforward manner
to any other bounds.)
\begin{proposition}\label{prop:bounded-square}
  Sceptic has a strategy in the bounded square-loss game
  that guarantees~(\ref{eq:level-3}).
\end{proposition}
\noindent
Proposition~\ref{prop:square} will be proved in Section~\ref{sec:level-2},
and it will imply Proposition~\ref{prop:bounded-square}
in combination with results of Section~\ref{sec:reductions}.

\subsection*{Counterexample}

The \emph{bounded absolute-loss game} is $([0,1],[0,1],\ell)$
where $\ell(\omega,\gamma):=\lvert\omega-\gamma\rvert$.
The level 3 Jeffreys theorem does not hold for the bounded absolute-loss game:
\begin{proposition}\label{prop:counter-3}
  Sceptic does not have a strategy
  that guarantees (\ref{eq:level-3}) in the bounded absolute-loss game.
\end{proposition}
\begin{proof}
  Suppose Sceptic has such a strategy and is playing it.
  Let Nature produce 0 and 1 independently with probability $1/2$ each.
  Predictor~1 always predicts $0$ and Predictor~2 always predicts $1$.
  The restriction of Sceptic's strategy
  to $\omega_n\in\{0,1\}$ and $\gamma_n^{[1]},\gamma_n^{[2]}\in\{0,1\}$
  is automatically measurable.
  We can see that $L_n^{[1]}-\tilde L_n$ and $L_n^{[2]}-\tilde L_n$
  are martingales with bounded increments,
  and so tend to $\infty$ with probability zero
  (see \cite{shiryaev:1996}, Theorem VII.5.1 and its corollary).
  Therefore, (\ref{eq:level-3}) happens with probability zero.
\end{proof}

The proof shows that Proposition~\ref{prop:counter-3}
remains true for the restricted game $(\{0,1\},[0,1],\ell)$,
$\ell(\omega,\gamma):=\lvert\omega-\gamma\rvert$.

\section{Reductions between Jeffreys theorems}
\label{sec:reductions}

It appears that the main factor that determines which Jeffreys theorems
hold for a particular game of prediction is the degree of convexity of the game.
We might define a game to be convex
if its prediction set $\Gamma$ is a convex set in a linear space
and its loss function $\ell(\omega,\gamma)$
is convex in $\gamma\in\Gamma$.
However, this definition would be too narrow,
since the predictions $\gamma$ are usually just arbitrary labels.
We start from introducing a much less arbitrary representation of games of prediction.

A \emph{canonical prediction} is a function $\lambda:\Omega\to\bbbr$
such that
\begin{equation*}
  \exists\gamma\in\Gamma
  \;
  \forall\omega\in\Omega:
  \lambda(\omega)
  =
  \ell(\omega,\gamma).
\end{equation*}
The \emph{canonical representation} of the game $(\Omega,\Gamma,\ell)$
is the pair $(\Omega,\Lambda)$
where $\Lambda$, called the \emph{canonical prediction set},
is the set of all canonical predictions.
We will not always distinguish between the game and its canonical representation
and will usually consider games that are non-redundant in the sense that
\begin{equation}\label{eq:non-redundant}
  (\lambda_1,\lambda_2\in\Lambda \;\&\; \lambda_1\le\lambda_2)
  \Longrightarrow
  \lambda_1=\lambda_2.
\end{equation}
A \emph{superprediction} (resp.\ \emph{subprediction})
is a function $\lambda:\Omega\to\bbbr$
such that $\lambda\ge\lambda'$ (resp.\ $\lambda\le\lambda'$)
for some canonical prediction $\lambda'$.
The set of all superpredictions (resp.\ subpredictions)
will be denoted $\upperLambda$ (resp.\ $\lowerLambda$)
and called the \emph{superprediction set} (resp.\ \emph{subprediction set}).

We will be interested in three notions of convexity for games of prediction:
\begin{itemize}
\item
  a game is \emph{convex} if its superprediction set $\upperLambda$ is convex
  (equivalently,
  if a convex mixture of two canonical predictions is always a superprediction);
  this condition is always satisfied
  if $\Gamma$ is a convex set
  and the loss function $\ell(\omega,\gamma)$ is convex in $\gamma\in\Gamma$;
\item
  a game is \emph{strictly convex}
  if a non-degenerate convex mixture of two canonical predictions
  is always an interior point of $\upperLambda$
  (in the topology of uniform convergence);
\item
  a game is \emph{perfectly mixable} if, for some $\eta>0$,
  the set $e^{-\eta\upperLambda}$ is convex.
\end{itemize}
For illustrative purposes
it is convenient to consider the case
where the game $(\Omega,\Gamma,\ell)$ is \emph{binary},
in the sense $\Omega=\{0,1\}$.
In this case $\Lambda$ can be represented as the subset of $\bbbr^2$
consisting of the points $(x,y)=(\lambda(0),\lambda(1))$
where $\lambda$ ranges over $\Lambda$.
An example is given as the curved line in Figure \ref{fig:hellinger} below;
the superpredictions are the points North-East of the line,
and the subpredictions are the points South-West of the line.


It is easy to see that for perfectly mixable prediction games
there is no real difference between the three levels of Jeffreys theorems:

\begin{proposition}\label{prop:equivalence}
  Suppose Sceptic can guarantee (\ref{eq:level-1})
  in the competitive prediction protocol
  for a perfectly mixable game.
  Then he can also guarantee (\ref{eq:level-3})
  (and, \emph{a fortiori}, (\ref{eq:level-2})).
\end{proposition}

\begin{proof}
  Consider the generalization of the competitive prediction protocol
  in which there are infinitely many Predictors
  (called Experts and numbered by $k=1,2,\ldots$)
  instead of just two.
  Using the Aggregating Algorithm
  (see, e.g., \cite{vovk:2001competitive}, Subsection~2.1),
  for any sequence $p_1,p_2,\ldots$
  of positive weights summing to $1$
  Sceptic can guarantee that his loss satisfies
  \begin{equation}\label{eq:AA}
    \tilde L_N
    \le
    L_N^{[k]}
    +
    C \ln\frac{1}{p_k}
  \end{equation}
  for all $N=1,2,\ldots$ and $k=1,2,\ldots$,
  where $C$ is a constant depending on the prediction game.

  Let Sceptic play a strategy that guarantees (\ref{eq:level-1}).
  We will construct a new strategy for Sceptic that guarantees (\ref{eq:level-3}).
  Consider the following doubly infinite set of experts:
  \begin{itemize}
  \item
    Expert $(k,1)$, $k=1,2,\ldots$, plays as Sceptic
    until the difference $L_n^{[1]}-\tilde L_n$ exceeds $2^k$;
    as soon as this happens (if it ever happens),
    he starts playing as Predictor 1;
  \item
    Expert $(k,2)$ plays as Sceptic
    until the difference $L_n^{[2]}-\tilde L_n$ exceeds $2^k$;
    as soon as this happens,
    he starts playing as Predictor 2.
  \end{itemize}
  The weights $p_{k,1}$ and $p_{k,2}$ assigned to these experts
  are $p_{k,1}=p_{k,2}=2^{-k-1}$.
  Applied to these experts,
  the Aggregating Algorithm provides a new strategy for Sceptic
  that guarantees (\ref{eq:level-3}).
  Indeed, suppose the first of the three terms in (\ref{eq:level-3})
  is false.
  Then, by (\ref{eq:level-1}),
  either the second or the third term in (\ref{eq:level-3}) becomes true
  when $\lim$ is replaced by $\limsup$.
  Suppose, for concreteness, it is the second term.
  For each $k$,
  Expert $(k,1)$'s loss satisfies
  $L_N^{[k,1]}<L_N^{[1]}-2^k$
  from some $N$ on,
  and so (\ref{eq:AA}) implies that the Aggregating Algorithm's loss $L_N$ satisfies
  $$
    L_N
    \le
    L_N^{[k,1]}
    +
    C \ln\frac{1}{p_{k,1}}
    <
    L_N^{[1]} - 2^k + (C\ln2)(k+1)
  $$
  for all $k$ and from some $N$ on.
  Letting $k\to\infty$,
  we can see that the second term of (\ref{eq:level-3}),
  with $L_N$ in place of $\tilde L_N$, is true.
\end{proof}

Of course,
Proposition~\ref{prop:equivalence} will continue to hold
if the Euclidean distance in (\ref{eq:level-1}), (\ref{eq:level-2}), and (\ref{eq:level-3})
is replaced by any other distance.

\subsection*{Examples of perfectly mixable games}

The bounded square-loss game is perfectly mixable
(\cite{vovk:2001competitive}, Subsection~2.4).

Perhaps the most fundamental class of games of prediction
is that of log-loss games.
If $(\Omega,\Gamma,\ell)$ is a \emph{log-loss game},
$\Omega$ is a measurable space with a fixed $\sigma$-finite measure $\mu$
(more generally, $\mu=\mu_n$ may depend on $n$
and be announced by a player, say Nature,
at the beginning of step $n$ of the game),
$\Gamma$ is the set of all measurable functions $\gamma:\Omega\to[0,\infty)$
satisfying $\int\gamma\dd\mu=1$,
and $\ell(\omega,\gamma)=-\ln\gamma(\omega)$.
For log-loss games the loss function is allowed to take value $\infty$
($-\ln0:=\infty$).
A simple and instructive special case to keep in mind
is where $\mu$ is the counting measure on a countable $\Omega$.
The perfect mixability of log-loss games is a well-known fact,
and the Aggregating Algorithm for them reduces to the Bayes rule
(for details see, e.g., \cite{vovk:2001competitive}, Subsection~2.2).

For other examples of perfectly mixable games
(such as the Kullback--Leibler game and Cover's game),
see \cite{vovk:2001competitive}, Subsection~2.5.

\section{Level 2 Jeffreys theorems}
\label{sec:level-2}

\begin{figure}[tb]
  \centering
    \includegraphics[width=10cm,clip=true]{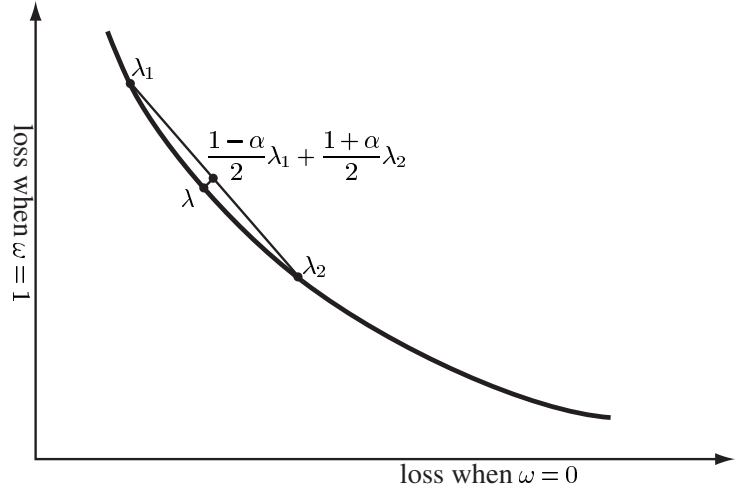}
  \caption{The interpretation of the $\alpha$-divergence
    between canonical predictions $\lambda_1$ and $\lambda_2$ in the binary case:
    find the mean $\frac{1-\alpha}{2}\lambda_1 + \frac{1+\alpha}{2}\lambda_2$
    of $\lambda_1$ and $\lambda_2$;
    find the intersection $\lambda$ of the prediction set
    and the slope $1$ line passing through the mean;
    multiply the horizontal (=vertical) distance between the mean and $\lambda$
    by $\frac{4}{1-\alpha^2}$.%
  \label{fig:hellinger}}
\end{figure}

If $\lambda_1$ and $\lambda_2$ are canonical predictions and $\alpha\in(-1,1)$,
we set
\begin{equation}\label{eq:lower-divergence}
  \lowerD^{[\alpha]}
  (\lambda_1\para\lambda_2)
  :=
  \frac{4}{1-\alpha^2}
  \sup
  \left\{
    t\in\bbbr
    \st
    \frac{1-\alpha}{2}
    \lambda_1
    +
    \frac{1+\alpha}{2}
    \lambda_2
    -
    t
    \in
    \upperLambda
  \right\}
\end{equation}
(the \emph{lower $\alpha$-divergence between $\lambda_1$ and $\lambda_2$})
and
\begin{equation*}
  \upperD^{[\alpha]}
  (\lambda_1\para\lambda_2)
  :=
  \frac{4}{1-\alpha^2}
  \inf
  \left\{
    t\in\bbbr
    \st
    \frac{1-\alpha}{2}
    \lambda_1
    +
    \frac{1+\alpha}{2}
    \lambda_2
    -
    t
    \in
    \lowerLambda
  \right\}
\end{equation*}
(the \emph{upper $\alpha$-divergence between $\lambda_1$ and $\lambda_2$}).
The lower and upper divergence make take values $-\infty$ or $\infty$.
We will be mostly interested in lower divergences
(which for many interesting games coincides with upper divergences).
In the case of binary $(\Omega,\Gamma,\ell)$
this definition is illustrated in Figure \ref{fig:hellinger}
(notice that the difference between lower and upper $\alpha$-divergences
disappears for convex binary games;
in such cases,
we will sometimes write $D^{[\alpha]}(\lambda_1\para\lambda_2)$
for the common value of $\lowerD^{[\alpha]}(\lambda_1\para\lambda_2)$
and $\upperD^{[\alpha]}(\lambda_1\para\lambda_2)$
and omit the adjectives ``lower'' and ``upper'').
We will also write $\lowerD^{[\alpha]}(\gamma_1\para\gamma_2)$
and $\upperD^{[\alpha]}(\gamma_1\para\gamma_2)$
for $\gamma_1,\gamma_2\in\Gamma$,
in the obvious sense.

Notice that, for strictly convex and non-redundant
(in the sense of (\ref{eq:non-redundant})) games,
\begin{equation*}
  \upperD^{[\alpha]}
  (\lambda_1\para\lambda_2)
  \ge
  \lowerD^{[\alpha]}
  (\lambda_1\para\lambda_2)
  >
  0,
\end{equation*}
for all $\lambda_1,\lambda_2\in\Lambda$.
For $\alpha=0$ the lower (resp.\ upper) $\alpha$-divergence
is called the \emph{lower} (resp.\ \emph{upper}) \emph{Hellinger distance};
the word ``distance'' is partly explained by its symmetry
(although simplest examples show that there is no continuous function $f$
such that $f(\lowerD^{[0]})$ or $f(\upperD^{[0]})$ is a metric
for every strictly convex game).

The values of lower and upper $\alpha$-divergences for $\alpha=\pm1$
are defined as their limits
as $\alpha\to\pm1$ when those limits exist.
The lower (resp.\ upper) $-1$-divergence is called
the \emph{lower} (resp.\ \emph{upper}) \emph{Kullback--Leibler divergence}
and is especially important.

\begin{remark}
  It is not difficult to see that upper divergences can be very different
  from the corresponding lower divergences
  even for ``nice'' (in particular, strictly convex) games.
  For example,
  for the game $([-1,1],[-1,1],(\omega-\gamma)^4)$
  the lower and upper Hellinger distances
  between the predictions $-1$ and $1$ are different, 1 and 7.
  (Cf.\ \cite{dawid:2007}, Lemma~3.)
\end{remark}

\subsection*{The square-loss and log-loss games}

In this subsections we will compute lower and upper divergences
for two popular games of prediction defined earlier.

\begin{lemma}\label{lem:divergence-square-loss}
  In the square-loss game,
  \begin{equation}\label{eq:divergence-square-loss}
    D^{[\alpha]}
    (\gamma_1\para\gamma_2)
    =
    (\gamma_1-\gamma_2)^2
  \end{equation}
  for all $\alpha\in[-1,1]$ and $\gamma_1,\gamma_2\in\bbbr$.
\end{lemma}
\begin{proof}
  It suffices to consider the case $\alpha\in(-1,1)$.
  The statement of the lemma will follow from the fact
  that, for all $\omega\in\bbbr$,
  \begin{multline*}
    \frac{1-\alpha}{2}
    (\gamma_1-\omega)^2
    +
    \frac{1+\alpha}{2}
    (\gamma_2-\omega)^2
    -
    \frac{1-\alpha^2}{4}
    (\gamma_1-\gamma_2)^2\\
    =
    \left(
      \frac{1-\alpha}{2}
      \gamma_1
      +
      \frac{1+\alpha}{2}
      \gamma_2
      -
      \omega
    \right)^2.
  \end{multline*}
  If we set $t_1:=\gamma_1-\omega$ and $t_2:=\gamma_2-\omega$,
  the last equality simplifies to the obvious
  \begin{equation*}
    \frac{1-\alpha}{2}
    t_1^2
    +
    \frac{1+\alpha}{2}
    t_2^2
    -
    \frac{1-\alpha^2}{4}
    (t_1-t_2)^2
    =
    \left(
      \frac{1-\alpha}{2}
      t_1
      +
      \frac{1+\alpha}{2}
      t_2
    \right)^2.
    \qedhere
  \end{equation*}
\end{proof}

\begin{lemma}\label{lem:divergence-log-loss}
  In any log-loss game,
  \begin{equation}\label{eq:divergence-log-loss}
    D^{[\alpha]}
    (\gamma_1\para\gamma_2)
    =
    -\frac{4}{1-\alpha^2}
    \ln
    \int_{\Omega}
      (\gamma_1(\omega))^{\frac{1-\alpha}{2}}
      (\gamma_2(\omega))^{\frac{1+\alpha}{2}}
    \mu(\dd\omega)
  \end{equation}
  for all $\alpha\in(-1,1)$ and $\gamma_1,\gamma_2\in\Gamma$.
\end{lemma}
\begin{proof}
  The left-hand side of (\ref{eq:divergence-log-loss})
  can be written as $\frac{4}{1-\alpha^2}t$
  where $t$ is defined from the condition
  that, for some $\gamma\in\Gamma$ and all $\omega\in\Omega$,
  \begin{equation*}
    -\frac{1-\alpha}{2}
    \ln\gamma_1(\omega)
    -
    \frac{1+\alpha}{2}
    \ln\gamma_2(\omega)
    -
    t
    =
    -\ln\gamma(\omega).
  \end{equation*}
  Deducing
  \begin{equation*}
  \int_{\Omega}\gamma\dd\mu
    =
    \int_{\Omega}
      (\gamma_1(\omega))^{\frac{1-\alpha}{2}}
      (\gamma_2(\omega))^{\frac{1+\alpha}{2}}
    \mu(\dd\omega)
    e^{t},
  \end{equation*}
  substituting $1$ for $\int\gamma\dd\mu$,
  and solving the resulting equation for $t$,
  we obtain the statement of the lemma.
\end{proof}

The standard definition of the $\alpha$-divergence
for the log-loss game
(see, e.g., \cite{amari/nagaoka:2000}, p.~57) is
\begin{equation*}
  D^{(\alpha)}
  (\gamma_1\para\gamma_2)
  =
  \frac{4}{1-\alpha^2}
  \left(
    1
    -
    \int_{\Omega}
      (\gamma_1(\omega))^{\frac{1-\alpha}{2}}
      (\gamma_2(\omega))^{\frac{1+\alpha}{2}}
    \mu(\dd\omega)
  \right);
\end{equation*}
it is clear that this will differ little from (\ref{eq:divergence-log-loss})
when $\gamma_1$ and $\gamma_2$ are close in a suitable sense.
The inequality $\ln x\le x-1$ implies $D^{(\alpha)}\le D^{[\alpha]}$.

\subsection*{Level 2 and level 3 Jeffreys theorems}

This is our most general level 2 Jeffreys theorem:

\begin{proposition}\label{prop:jeffreys}
  For each $\alpha\in(-1,1)$ and $\epsilon>0$
  Sceptic has a strategy that guarantees
  \begin{equation}\label{eq:jeffreys}
    \frac{1-\alpha^2}{4}
    \sum_{n=1}^N
    \lowerD^{[\alpha]}
    \left(
      \gamma^{[1]}_n\para\gamma^{[2]}_n
    \right)
    \le
    \frac{1-\alpha}{2}
    L^{[1]}_N
    +
    \frac{1+\alpha}{2}
    L^{[2]}_N
    -
    \tilde L_N
    +
    \epsilon
  \end{equation}
\end{proposition}
\begin{proof}
  The strategy is obvious:
  according to (\ref{eq:lower-divergence}),
  at step $n$ Sceptic can choose a canonical prediction $\lambda$ satisfying
  \begin{equation*}
    \lambda
    \le
    \frac{1-\alpha}{2}
    \lambda_1
    +
    \frac{1+\alpha}{2}
    \lambda_2
    -
    \frac{1-\alpha^2}{4}
    \lowerD^{[\alpha]}
    (\lambda_1\para\lambda_2)
    +
    \epsilon 2^{-n}
  \end{equation*}
  ($\lambda_1$ and $\lambda_2$ being the canonical predictions
  corresponding to $\gamma_n^{[1]}$ and $\gamma_n^{[2]}$).
  Summing over the first $N$ steps,
  we obtain (\ref{eq:jeffreys}).
\end{proof}

Specializing (\ref{eq:jeffreys}) to the case $\alpha=0$
and the square-loss game gives
\begin{equation*}
  \frac{1}{4}
  \sum_{n=1}^N
  \left(
    \gamma^{[1]}_n-\gamma^{[2]}_n
  \right)^2
  \le
  \frac
  {
    L^{[1]}_N
    +
    L^{[2]}_N
  }
  {2}
  -
  \tilde L_N
  +
  \epsilon.
\end{equation*}
This implies a stronger version of the level 2 Jeffreys theorem (\ref{eq:level-2}):
\begin{equation*} 
  \sum_{n=1}^{\infty}
  \left(
    \gamma^{[1]}_n-\gamma^{[2]}_n
  \right)^2
  <
  \infty
  \text{ or }
  \lim_{N\to\infty}
  \max
  \left(
    L^{[1]}_N
    -
    \tilde L_N,
    L^{[2]}_N
    -
    \tilde L_N
  \right)
  =
  \infty.
\end{equation*}
In combination with the proof of Proposition~\ref{prop:equivalence},
this implies the stronger form
\begin{equation}\label{eq:level-3-strong}
  \sum_{n=1}^{\infty}
  \left(
    \gamma^{[1]}_n-\gamma^{[2]}_n
  \right)^2
  <
  \infty
  \text{ or }
  \lim_{N\to\infty}
  \left(
    L^{[1]}_N
    -
    \tilde L_N
  \right)
  =
  \infty
  \text{ or }
  \lim_{N\to\infty}
  \left(
    L^{[2]}_N
    -
    \tilde L_N
  \right)
  =
  \infty
\end{equation}
of the level 3 Jeffreys theorem (\ref{eq:level-3})
for the bounded square-loss game.

For the log-loss game,
we obtain (\ref{eq:level-3-strong})
with the Hellinger distance
$D^{[0]}(\gamma^{[1]}_n\para\gamma^{[2]}_n)$,
or the standard Hellinger distance
$D^{(0)}(\gamma^{[1]}_n\para\gamma^{[2]}_n)$,
in place of $(\gamma^{[1]}_n-\gamma^{[2]}_n)^2$.

\section{Level 1 Jeffreys theorems}
\label{sec:level-1}

The main goal of this section is to prove Proposition~\ref{prop:absolute}.
In the absolute-loss game,
the divergence between any two predictions is $0$,
and so the methods of the previous section are not applicable.

\begin{figure}[tb]
    \centering
      \includegraphics[width=10cm,clip=true]{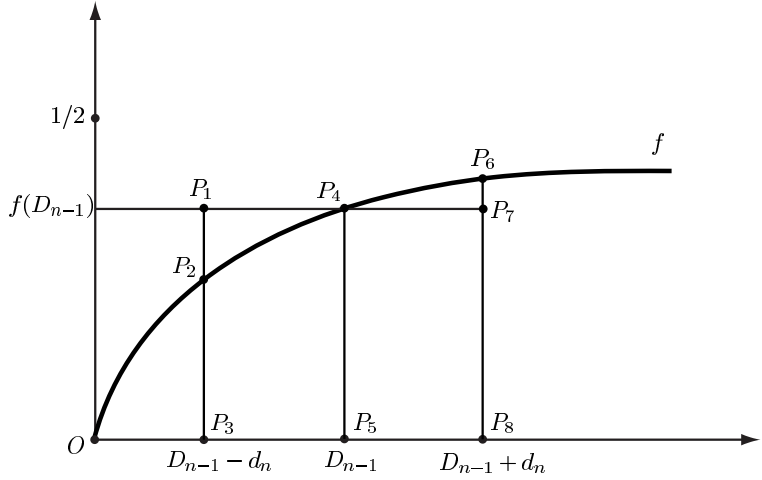}
    \caption{The function $f$ from the proof of Proposition~\ref{prop:absolute}.%
    \label{fig:absolute}}
  \end{figure}

  First we describe a strategy for Sceptic
  that will later be shown to ensure (\ref{eq:level-1}).
  Let $f:[0,\infty)\to[0,1/2)$ be a strictly increasing and concave function
  satisfying $f(0)=0$ and $f(\infty)<1/2$;
  see Figure \ref{fig:absolute}.
  Later it will be convenient to extend $f$ to $(-\infty,\infty)$
  by the central symmetry w.r.\ to the origin $O$
  (so that $f:(-\infty,\infty)\to(-1/2,1/2)$ is an odd function).

  Suppose just before step $n=1,2,\ldots$ of the competitive prediction protocol
  we have $D_{n-1}:=L_{n-1}^{[1]}-L_{n-1}^{[2]}\ge0$
  (the case where $L_{n-1}^{[1]}\le L_{n-1}^{[2]}$
  will later be reduced to this one).
  Sceptic's move can be represented as
  $$
    \tilde\gamma_n
    :=
    (1-t_n) \gamma_n^{[1]}
    +
    t_n \gamma_n^{[2]},
  $$
  where $t_n$ will be chosen later from the interval $[0,1/2]$.
  Set
  \begin{align*}
    d_n &:= \left|\gamma_n^{[1]} - \gamma_n^{[2]}\right| \in [0,1],\\
    \bar L_n &:= \frac{L_n^{[1]} + L_n^{[2]}}{2},\\
    \bar \ell_n &:= \frac{\ell(\omega_n,\gamma_n^{[1]}) + \ell(\omega_n,\gamma_n^{[2]})}{2}.
  \end{align*}

  If the actual outcome $\omega_n$ is in favour of Predictor 1,
  $$
    \ell(\omega_n,\gamma_n^{[1]})
    \le
    \ell(\omega_n,\gamma_n^{[2]}),
  $$
  the difference $L_n^{[1]} - L_n^{[2]}$ between the losses of the two Predictors
  will decrease to $D_n=D_{n-1}-d_n$
  and the difference $\tilde L_n - \bar L_n$ will increase by
  \begin{equation*}
    \ell(\omega_n,\tilde\gamma_n)
    -
    \bar\ell_n
    =
    (1-t_n)
    \left(
      \bar\ell_n - \frac{d_n}{2}
    \right)
    +
    t_n
    \left(
      \bar\ell_n + \frac{d_n}{2}
    \right)
    -
    \bar\ell_n
    =
    \left(
      t_n - \frac{1}{2}
    \right)
    d_n.
  \end{equation*}
  So in fact it will decrease as $t_n\le1/2$.
  Let us set $t_n:=1/2-f(D_{n-1})$.
  The difference $\tilde L_n-\bar L_n$ will decrease
  by the area of the rectangle $P_3 P_5 P_4 P_1$.

  If the actual outcome $\omega_n$ is in favour of Predictor 2,
  $$
    \ell(\omega_n,\gamma_n^{[1]})
    \ge
    \ell(\omega_n,\gamma_n^{[2]}),
  $$
  the difference between the losses of the two Predictors
  will increase to $D_n=D_{n-1}+d_n$
  and the difference $\tilde L_n - \bar L_n$ will increase by
  \begin{multline*}
    \ell(\omega_n,\tilde\gamma_n)
    -
    \bar\ell_n
    =
    (1-t_n)
    \left(
      \bar\ell_n + \frac{d_n}{2}
    \right)
    +
    t_n
    \left(
      \bar\ell_n - \frac{d_n}{2}
    \right)
    -
    \bar\ell_n\\
    =
    \left(
      \frac{1}{2} - t_n
    \right)
    d_n
    =
    f(D_{n-1}) d_n,
  \end{multline*}
  i.e., by the area of the rectangle $P_5 P_8 P_7 P_4$.

  We can see that in both cases,
  $D_n=D_{n-1}\pm d_n$,
  the difference $\tilde L_n-\bar L_n$
  increases by $\int_{D_{n-1}}^{D_n}f$
  minus the area $A_n$ of a curvilinear triangle
  ($P_1P_2P_4$ if $D_n=D_{n-1}-d_n$
  and $P_4P_7P_6$ if $D_n=D_{n-1}+d_n$).
  Now extend $f$ to the whole of $(-\infty,\infty)$ as an odd function.
  Suppose that $D_{n-1}\le0$ and, moreover, $D_{n-1}+d_n\le0$.
  Applying the same argument as above
  but with the roles of Predictor 1 and Predictor 2 interchanged,
  we can see that the difference $\tilde L_n-\bar L_n$
  again increases by $\int_{D_{n-1}}^{D_n}f$
  minus the area $A_n$ of a curvilinear triangle.
  It is easy to check that the difference $\tilde L_n-\bar L_n$
  will change in the same way also
  in the case where $D_{n-1}\ge0$ but $D_{n-1}-d_n\le0$
  and in the case where $D_{n-1}\le0$ but $D_{n-1}+d_n\ge0$.
  Since $\tilde L_N-\bar L_N$ is the cumulative increase in $\tilde L_n-\bar L_n$
  over $n=1,\ldots,N$,
  we can see that
  $$
    \tilde L_N-\bar L_N
    =
    \int_0^{D_N}f
    -
    \sum_{n=1}^{N}
    A_n.
  $$
  It remains to consider two cases:
  \begin{description}
  \item[$\sum_{n=1}^\infty A_n<\infty$:]
    In this case,
    $A_N\to0$ and so
    $$
      \max
      \left(
        \frac{1}{\left|\gamma_N^{[1]}-\gamma_N^{[2]}\right|},
        \left|D_N\right|
      \right)
      \to
      \infty
    $$
    as $N\to\infty$.
    The sequence $N=1,2,\ldots$ can be split into three subsequences
    such that
    $\lvert\gamma_N^{[1]}-\gamma_N^{[2]}\rvert\to0$ along the first,
    $D_N\to\infty$ along the second,
    and $D_N\to-\infty$ along the third.
    It suffices to show that (\ref{eq:level-1})
    holds along the second subsequence
    (the case of the third subsequence is analogous,
    and the case of the first subsequence is trivial).
    Assuming $D_N>0$,
    we can see that along the second subsequence:
    \begin{multline*}
      \tilde L_N
      =
      \bar L_N
      +
      \int_0^{D_N}f
      -
      \sum_{n=1}^{N}
      A_n\\
      \le
      \frac{L_N^{[1]}+L_N^{[1]}-D_N}{2}
      +
      \int_0^{D_N}f
      \le
      L_N^{[1]}
      +
      D_N
      \left(
        f(\infty) - \frac12
      \right),
    \end{multline*}
    and so $L_N^{[1]}-\tilde L_N\to\infty$.
  \item[$\sum_{n=1}^\infty A_n=\infty$:]
    In this case we have along the subsequence of $N$ for which $D_N\ge0$:
    \begin{multline*}
      \tilde L_N
      =
      \bar L_N
      +
      \int_0^{D_N}f
      -
      \sum_{n=1}^{N}
      A_n\\
      =
      \frac{L_N^{[1]}+L_N^{[1]}-D_N}{2}
      +
      \int_0^{D_N}f
      -
      \sum_{n=1}^{N}
      A_n
      \le
      L_N^{[1]}
      -
      \sum_{n=1}^{N}
      A_n,
    \end{multline*}
    and so $L_N^{[1]}-\tilde L_N\to\infty$.
    Similarly,
    $L_N^{[2]}-\tilde L_N\to\infty$
    along the subsequence of $N$ for which $D_N\le0$.
    Therefore, (\ref{eq:level-1}) holds.
  \end{description}

\subsection*{Convex games}

It is easy to see that the proof of Proposition~\ref{prop:absolute}
is applicable to any convex game.
For any such game
Sceptic has a strategy in the competitive prediction protocol
that guarantees
\begin{equation*}
  \lim_{N\to\infty}
  \max
  \left(
    \frac{1}
    {
      \left|
        \lambda\left(\omega_N,\gamma^{[1]}_N\right)
        -
        \lambda\left(\omega_N,\gamma^{[2]}_N\right)
      \right|
    },
    L^{[1]}_N
    -
    \tilde L_N,
    L^{[2]}_N
    -
    \tilde L_N
  \right)
  =
  \infty.
\end{equation*}

\subsection*{Acknowledgements}

We are grateful to Akio Fujiwara for a useful discussion,
to Glenn Shafer for his advice,
and to participants of WITMSE 2009 for their comments.
This work was supported in part by EPSRC (grant EP/F002998/1).

\end{document}